%% file: Krein.tex
\documentclass[12pt]{article} 
\usepackage{amsfonts}          
\usepackage{myart2k}
\title{\ \\ \normalsize\bf{ON $J$-CONSERVATIVE SCATTERING SYSTEM
 REALIZATION \\ IN SEVERAL VARIABLES}}
\author{\normalsize\bf{D. S. Kalyuzhniy-Verbovetzky}}
\date{}
\newcommand{\nspace}[2]{\ensuremath{{\mathbb{#1}}^{#2}}}
\newcommand{\Hspace}[1]{\ensuremath{\mathcal{#1}}}
\input{mydef}

\begin{document}
\maketitle \vspace{-1cm}
\begin{abstract}
\noindent We prove that an arbitrary function, which is
holomorphic on some neighbourhood of $z=0$ in $\mathbb{C}^N$ and
vanishes at $z=0$, whose values are bounded linear operators
mapping one separable Hilbert space into another one, can be
represented as the transfer function of some multiparametric
discrete time-invariant conservative scattering linear system with
a Krein space of its inner states.
\end{abstract}
\thispagestyle{empty} \setcounter{section}{-1}
\section{Introduction}\label{sec:intr}
An arbitrary function, holomorphic on a neighbourhood of $z=0$ in
$\mathbb{C}$, whose values are bounded linear operators mapping
one separable Hilbert space into another one, can be represented
as the transfer fuction of some discrete time-invariant linear
system. This fact was established by D.Z.~Arov in \cite{A_1974}
(see also \cite{A1_1979,BGK_1979}). Later on, T.Ya.~Azizov proved
(in a different terminology) that an arbitrary holomorphic
operator-valued function on a neighbourhood of $z=0$ in
$\mathbb{C}$ has a $J$-conservative scattering system realization
(see \cite{AI_1986}). For that, he constructed a $J$-conservative
scattering system dilation of an arbitrary discrete time-invariant
linear system. Then, using realization mentioned earlier and
taking its $J$-conservative scattering system dilation, due to the
fact that a transfer function remains the same under dilation of a
system, he obtained a desired $J$-conservative realization.

Let us note that there exist other proofs of this $J$-conservative
realization theorem (see, e.g., \cite{ADRdS_1997} and references
there), however in this paper we shall follow, in the main, the
same scheme as in Azizov's proof, for the proof of the
multivariable generalization of this theorem.

First, we introduce the notion of multiparametric discrete
time-invariant $J$-conser\-vative scattering linear system (or
conservative scattering linear system with a Krein space of inner
states) which generalizes the notion of multiparametric
conservative scattering linear system introduced in
\cite{K1_2000}. Second, we prove that an arbitrary multiparametric
linear system has a $J$-conservative scattering system dilation
(the notion of dilation for multiparametric linear systems was
introduced in \cite{K2_2000}). Next, we use our generalization of
the realization theorem mentioned in the very beginning of this
paper, to several variables \cite{K3_2000}, and the fact that the
transfer function of a multiparametric linear system remains the
same under dilation \cite{K2_2000}, and prove the main result of
this paper on the existence of a multiparametric $J$-conservative
scattering system realization for an arbitrary operator-valued
function which is holomorphic on some neighbourhood of $z=0$ in
$\nspace{C}{N}$ and vanishes at $z=0$.

The organization of this paper is as follows. In
Section~\ref{sec:prelim} we give some preliminaries, with exact
formulations of the results mentioned above, on $J$-conservative
scattering linear systems in the one-parametric discrete case, in
the system-theoretical language convenient for the sequel. In
Section~\ref{sec:mult} we introduce multiparametric
$J$-conservative scattering linear systems (in the discrete case),
and formulate our main theorems. Section~\ref{sec:proofs} contains
the proofs of these theorems.

\section{Preliminaries on one-parametric $J$-conservative scattering
 linear systems}\label{sec:prelim}
 In our notation, a \emph{one-parametric discrete time-invariant
 linear system $\alpha$} is given by
\begin{equation}\label{eq:1-sys}
\alpha :\left\{ \begin{array}{lll} x(t)&=&Ax(t-1)+Bu(t-1),\\
y(t)&=&Cx(t-1)+Du(t-1)
\end{array}
\right.\quad (t=1,2,\ldots ),
\end{equation}
 where for each $t$ vectors $x(t),u(t),y(t)$ belong to separable
 Hilbert spaces $\Hspace{X}, \Hspace{U},
\Hspace{Y}$, respectively (throughout this paper we consider only
such type of spaces when nothing is said especially);
$A:\Hspace{X}\to\Hspace{X},\ B:\Hspace{U}\to\Hspace{X},\
C:\Hspace{X}\to\Hspace{Y},\ D:\Hspace{U}\to\Hspace{Y}$ are bounded
linear operators.

Such a form of a system differs from the standard one by the unit
shift in the argument of an output signal $y(\cdot)$, that leads,
in fact, to the equivalent theory (for more details and motivation
of such a notation of a system, see \cite{K1_2000}). Thus, for
example, the \emph{transfer function} of a system $\alpha$ of the
form \eqref{eq:1-sys} is given by
\begin{equation}\label{eq:1-tf}
\theta_\alpha(z)=zD+zC(I_{\Hspace{X}}-zA)^{-1}zB
\end{equation}
in some neighbourhood of $z=0$ in $\mathbb{C}$, i.e. differs from
the standard one by multiplier $z$ only (here $I_\Hspace{X}$
denotes the identity operator on $\Hspace{X}$).

The first result mentioned in Section~\ref{sec:intr} can be
formulated now as follows: an arbitrary
$L(\Hspace{U},\Hspace{Y})$-valued function $\theta$ holomorphic on
some neighbourhhod $\Gamma$ of $z=0$ in $\mathbb{C}$ and vanishing
at $z=0$ can be realized as the transfer function of some system
$\alpha$ of the form \eqref{eq:1-sys}, i.e.,
$\theta(z)=\theta_\alpha (z)$ in some neighbourhood (possibly,
smaller than $\Gamma$) of $z=0$ (here we denote by
$L(\Hspace{U},\Hspace{Y})$ the Banach space of all bounded linear
operators from a separable Hilbert space $\Hspace{U}$ into a
separable Hilbert space $\Hspace{Y}$).

Let the operator $J\in L(\Hspace{X}):=L(\Hspace{X},\Hspace{X})$ be
given such that $J=J^*=J^{-1}$ (such a $J$ is said to be a
\emph{canonical symmetry on $\Hspace{X}$}). Then $J$ determines on
$\Hspace{X}$ the new inner product $[x_1,x_2]_J:=\langle
Jx_1,x_2\rangle$ (here $\langle\cdot,\cdot\rangle$ stands for a
Hilbert space inner product on $\Hspace{X}$) which is, in general,
indefinite, and the space $\Hspace{X}$ with this new inner product
has the structure of a Krein space (for more information on Krein
spaces see, e.g., \cite{AI_1986}).

Let $\alpha=(1;A,B,C,D;\Hspace{X},\Hspace{U},\Hspace{Y})$ be a
one-parametric linear system of the form \eqref{eq:1-sys}, and
$J\in L(\Hspace{X})$ be a canonical symmetry. Set $J_1:=J\oplus
I_{\Hspace{U}}\in L(\Hspace{X}\oplus\Hspace{U}),\ J_2:=J\oplus
I_{\Hspace{Y}}\in L(\Hspace{X}\oplus\Hspace{Y})$. We shall call
$\alpha$ a \emph{one-parametric $J$-conservative scattering
system} if the system operator
\begin{displaymath}
G=\left(\begin{array}{cc}
  A & B \\
  C & D
\end{array}\right)\in L(\Hspace{X}\oplus\Hspace{U},\Hspace{X}
\oplus\Hspace{Y})
\end{displaymath}
is ($J_1,J_2$)-unitary, i.e.
\begin{displaymath}
 G^*J_2G=J_1,\quad GJ_1G^*=J_2.
\end{displaymath}
In the particular case when $J=I_{\Hspace{X}}$, a $J$-conservative
scattering system is a conservative scattering one.

Let us note that one may consider a $J$-conservative scattering
system $\alpha=(1;A,B,\- C,D;\Hspace{X},\Hspace{U},\Hspace{Y})$ as
a conservative scattering one, however with a Krein space of its
inner states, i.e., the equations of energy balance for $\alpha$
have the same form as for conservative scattering system with a
Hilbert space of inner states, but with $J$-metric
$[\cdot,\cdot]_J$ instead of Hilbert metric
$\langle\cdot,\cdot\rangle$ for states of a system $\alpha$:
\begin{displaymath}
[x(t),x(t)]_J-[x(t-1),x(t-1)]_J=\| u(t-1)\|^2-\| y(t)\|^2\quad
(t=1,2,\ldots),
\end{displaymath}
and the analogous equation holds for states, inputs and outputs of
the conjugate system
$\alpha^*:=(1;A^*,C^*,B^*,D^*;\Hspace{X},\Hspace{Y},\Hspace{U})$.

Recall \cite{A_1974} (see also \cite{A1_1979,BGK_1979}) that a
system $\widetilde{\alpha}
=(1;\widetilde{A},\widetilde{B},\widetilde{C},D;\widetilde{\Hspace{X}},
\Hspace{U}, \Hspace{Y})$ is called  a \emph{dilation of a system}
$\alpha =(1;A, B, C, D; \Hspace{X}, \Hspace{U}, \Hspace{Y})$ if
$\widetilde{\Hspace{X}}\supset\Hspace{X}$, and there exist
subspaces $\Hspace{D}$ and $\Hspace{D}_*$ in
$\widetilde{\Hspace{X}}$ such that
\begin{eqnarray*}
\widetilde{\Hspace{X}}=\Hspace{D}\oplus\Hspace{X}\oplus\Hspace{D_*},
\quad \widetilde{A}\Hspace{D}\subset\Hspace{D}, &
\widetilde{C}\Hspace{D}=\{ 0\} , &
{\widetilde{A}}^*\Hspace{D_*}\subset\Hspace{D_*},\quad
\widetilde{B}^*\Hspace{D_*}=\{ 0\},\\
A=P_\Hspace{X}\widetilde{A}|\Hspace{X}, &
B=P_\Hspace{X}\widetilde{B}, & C=\widetilde{C}|\Hspace{X}
\end{eqnarray*}
(here $P_\Hspace{X}$ stands for the orthogonal projector onto
$\Hspace{X}$ in $\widetilde{\Hspace{X}}$). For that, the transfer
functions of $\alpha$ and $\widetilde{\alpha}$ coincide in some
neighbourhood of $z=0$ in $\mathbb{C}$.

Now Azizov's result mentioned in Section~\ref{sec:intr} can be
formulated in the following way. An arbitrary system
$\alpha=(1;A,B,C,D;\Hspace{X},\Hspace{U},\Hspace{Y})$ has a
dilation $\widetilde{\alpha}
=(1;\widetilde{A},\widetilde{B},\widetilde{C},D;
\widetilde{\Hspace{X}}, \Hspace{U}, \Hspace{Y})$ which is a
one-parametric $J$-conservative scattering system for certain
canonical symmetry $J\in L(\widetilde{\Hspace{X}})$. As a
corollary (see Section~\ref{sec:intr}), an arbitrary
$L(\Hspace{U},\Hspace{Y})$-valued function $\theta$ holomorphic on
some neighbourhood $\Gamma$ of $z=0$ in $\mathbb{C}$ and vanishing
at $z=0$ can be realized as the transfer function of some system
$\alpha=(1;A,B,C,D;\Hspace{X},\Hspace{U},\Hspace{Y})$ of the form
\eqref{eq:1-sys}, which is a one-parametric $J$-conservative
scattering system for certain canonical symmetry $J\in
L(\Hspace{X})$, i.e., $\theta(z)=\theta_\alpha (z)$ in some
neighbourhood (possibly, smaller than $\Gamma$) of $z=0$.

\section{Multiparametric $J$-conservative scattering linear systems}
\label{sec:mult} Let us recall some definitions from
\cite{K1_2000} concerning \emph{multiparametric discrete
time-invariant linear systems}. Such a system $\alpha$ is given by
\begin{equation}\label{eq:n-sys}
\alpha :\left\{\begin{array}{r} x(t)=\sum\limits
_{k=1}^N(A_kx(t-e_k)+B_ku(t-e_k)),\\ y(t)=\sum\limits
_{k=1}^N(C_kx(t-e_k)+D_ku(t-e_k))
\end{array}
\right. \quad (t\in\nspace{Z}{N},\ |t|>0),
\end{equation}
where for $t=(t_1,\ldots,t_N)\in\nspace{Z}{N}$ we set
$|t|:=\sum_{k=1}^Nt_k$, for all $k\in\{ 1,\ldots ,N\}$ we set
$e_k:=(0,\ldots,0,1,0,\ldots,0)\in \nspace{Z}{N}$ with $1$ on the
$k$-th place, and zeros on other places; for each $t$ vectors
$x(t),u(t),y(t)$ belong to (separable Hilbert) spaces
 $\Hspace{X}, \Hspace{U}, \Hspace{Y}$, respectively;
  for all $k\in\{1,\ldots ,N\}\
A_k:\Hspace{X}\to\Hspace{X},\ B_k:\Hspace{U}\to\Hspace{X},\
C_k:\Hspace{X}\to\Hspace{Y},\ D_k:\Hspace{U}\to\Hspace{Y}$ are
bounded linear operators. We shall use the short notation $\alpha
=(N; \mathbf{A}, \mathbf{B}, \mathbf{C}, \mathbf{D};\Hspace{X},
\Hspace{U}, \Hspace{Y})$ where $\mathbf{A}, \mathbf{B},
\mathbf{C}, \mathbf{D}$ mean $N$-tuples of operators
$A_k,B_k,C_k,D_k$, respectively. The \emph{transfer function} of a
system $\alpha$ of the form \eqref{eq:n-sys} is given by
\begin{equation}\label{eq:n-tf}
\theta _\alpha
(z)=z\mathbf{D}+z\mathbf{C}{(I_\Hspace{X}-z\mathbf{A})}^{-1}z
\mathbf{B}
\end{equation}
in some neighbourhood of $z=0$ in $\nspace{C}{N}$, where for
$z=(z_1,\ldots ,z_N)\in\nspace{C}{N}$ and an $N$-tuple of
operators $\mathbf{T}=(T_1,\ldots,T_N)$ we use the notation
$z\mathbf{T}:=\sum_{k=1}^Nz_kT_k$. It is clear that system
\eqref{eq:n-sys} is the generalization of system \eqref{eq:1-sys},
and formula \eqref{eq:n-tf} is the generalization of formula
\eqref{eq:1-tf} for transfer function, to the case of several
variables.

Recall \cite{K2_2000} that a system $\widetilde{\alpha}
=(N;\widetilde{\mathbf{A}}, \widetilde{\mathbf{B}},
\widetilde{\mathbf{C}}, \mathbf{D};\widetilde{\Hspace{X}},
\Hspace{U}, \Hspace{Y})$ is called a \emph{dilation of a
multiparametric linear system} $\alpha =(N; \mathbf{A},
\mathbf{B}, \mathbf{C}, \mathbf{D};\Hspace{X}, \Hspace{U},
\Hspace{Y})$ if for each $z\in\nspace{C}{N}$ a system
$\widetilde{\alpha}_z:=(1; z\widetilde{\mathbf{A}},
z\widetilde{\mathbf{B}}, z\widetilde{\mathbf{C}},
z\mathbf{D};\widetilde{\Hspace{X}}, \Hspace{U}, \Hspace{Y})$ is a
dilation of a one-parametric linear system $\alpha_z:=(1;
z\mathbf{A}, z\mathbf{B}, z\mathbf{C}, z\mathbf{D};\Hspace{X},
\Hspace{U}, \Hspace{Y})$, i.e.,
$\widetilde{\Hspace{X}}\supset\Hspace{X}$, and there exist
subspaces $\Hspace{D}_z$ and $\Hspace{D}_{*,z}$ in
$\widetilde{\Hspace{X}}$ such that
\begin{eqnarray*}
 \widetilde{\Hspace{X}}=\Hspace{D}_z\oplus\Hspace{X}\oplus
\Hspace{D}_{*,z},\quad
z\widetilde{\mathbf{A}}\Hspace{D}_z\subset\Hspace{D}_z, &
z\widetilde{\mathbf{C}}\Hspace{D}_z =\{ 0\}, &
(z\widetilde{\mathbf{A}})^*\Hspace{D}_{*,z}\subset
\Hspace{D}_{*,z},\quad
(z\widetilde{\mathbf{B}})^*\Hspace{D}_{*,z}=\{ 0\} \\
z\mathbf{A}=P_\Hspace{X}(z\widetilde{\mathbf{A}})|\Hspace{X}, &
z\mathbf{B}=P_\Hspace{X}(z\widetilde{\mathbf{B}}), &
z\mathbf{C}=(z\widetilde{\mathbf{C}})|\Hspace{X}.
\end{eqnarray*}
For that, the transfer functions of $\alpha$ and
$\widetilde{\alpha}$ coincide in some neighbourhood of $z=0$ in
$\nspace{C}{N}$.

Let  $\nspace{T}{N}:=\{\zeta\in\nspace{C}{N}:|\zeta_k|=1,\
k=1,\ldots,N\}$ be the $N$-dimensional unit torus.
\begin{defn}\label{defn:j-cons}
Let $\alpha =(N; \mathbf{A}, \mathbf{B}, \mathbf{C},
\mathbf{D};\Hspace{X}, \Hspace{U}, \Hspace{Y})$ and a canonical
symmetry $J\in L(\Hspace{X})$ be given. We shall call $\alpha$ a
\emph{multiparametric $J$-conservative scattering linear system}
if for each $\zeta\in\nspace{T}{N}$ $\alpha_\zeta:=(1;
\zeta\mathbf{A}, \zeta\mathbf{B}, \zeta\mathbf{C},
\zeta\mathbf{D};\Hspace{X}, \Hspace{U}, \Hspace{Y})$ is a
one-parametric $J$-conservative scattering linear system, i.e. for
$J_1=J\oplus I_{\Hspace{U}}\in L(\Hspace{X}\oplus\Hspace{U}),\
J_2=J\oplus I_{\Hspace{Y}}\in L(\Hspace{X}\oplus\Hspace{Y})$ one
has
\begin{equation}\label{eq:n-j-cons}
  (\zeta\mathbf{G})^*J_2(\zeta\mathbf{G})=J_1,\quad
  (\zeta\mathbf{G})J_1(\zeta\mathbf{G})^*=J_2,
\end{equation}
where $\mathbf{G}=(G_1,\ldots,G_N)$ is the $N$-tuple of system
operators
\begin{equation}\label{eq:gk}
  G_k=\left(\begin{array}{cc}
    A_k & B_k \\
    C_k & D_k
  \end{array}\right):\Hspace{X}\oplus\Hspace{U}\to
  \Hspace{X}\oplus\Hspace{Y} \quad (k=1,\ldots,N).
\end{equation}
\end{defn}

In the particular case when $J=I_{\Hspace{X}}$, this notion
coincides with the notion of multiparametric conservative
scattering linear system introduced in \cite{K1_2000}. Let us
remark here that another type of multiparametric conservative
scattering linear systems for the discrete case was considered by
J.A.~Ball and T.T.~Trent in \cite{BT_1998}.

By equating corresponding coefficients of trigonometric
polynomials in the left and right parts of \eqref{eq:n-j-cons},
one can easily see that $\alpha =(N; \mathbf{A}, \mathbf{B},
\mathbf{C}, \mathbf{D};\Hspace{X}, \Hspace{U}, \Hspace{Y})$ is a
multiparametric $J$-conservative scattering linear system if and
only if the following equalities hold:
\begin{eqnarray}
  \sum\limits_{k=1}^NG_k^*J_2G_k &=& J_1, \label{eq:j-c1} \\
G_k^*J_2G_l &=& 0 \quad (k\neq l), \label{eq:j-c2} \\
\sum\limits_{k=1}^NG_kJ_1G_k^* &=& J_2, \label{eq:j-c3} \\
G_kJ_1G_l^* &=& 0 \quad (k\neq l). \label{eq:j-c4}
\end{eqnarray}

This definition can be also formulated in terms of energy balance
equations, i.e. the conservation of energy for a system $\alpha
=(N;\mathbf{A}, \mathbf{B}, \mathbf{C}, \mathbf{D};\Hspace{X},
\Hspace{U}, \Hspace{Y})$, and for its \emph{conjugate system}
$\alpha^*:=(N; \mathbf{A}^*, \mathbf{C}^*, \mathbf{B}^*,
\mathbf{D}^*;\Hspace{X}, \Hspace{Y}, \Hspace{U})$, where for
$N$-tuples $\mathbf{T}=(T_1,\ldots,T_N)$ of operators $T_k\in
L(\Hspace{H}_1,\Hspace{H}_2)$ we set
$\mathbf{T}^*:=(T_1^*,\ldots,T_N^*)$ with $T_k^*\in
L(\Hspace{H}_2,\Hspace{H}_1),\ k=1,\ldots,N$, and for the "energy"
of states of $\alpha$ and $\alpha^*$ use $J$-metric
$[\cdot,\cdot]_J$ instead of Hilbert metric
$\langle\cdot,\cdot\rangle$. More precisely, the following
proposition is valid.
\begin{prop} $\alpha
=(N;\mathbf{A}, \mathbf{B}, \mathbf{C}, \mathbf{D};\Hspace{X},
\Hspace{U}, \Hspace{Y})$ is a multiparametric $J$-conservative
scattering linear system with some canonical symmetry $J\in
L(\Hspace{X})$ if and only if
\begin{description}
  \item[(i)] for any input multisequence $\{ u(t):|t|\geq 0\}$ of
  $\alpha$ satisfying the condition
  \begin{displaymath}
  \forall n\in\mathbb{N}\quad \sum\limits_{|t|=n}\| u(t)\|^2<\infty,
  \end{displaymath}
  and its initial states collection $\{ x(t):|t|=0\}$ satisfying
  the condition
\begin{displaymath}
   \sum\limits_{|t|=0}\| x(t)\|^2<\infty,
  \end{displaymath}
  one has for any $n\in\mathbb{N}: \
  \sum_{|t|=n}\| x(t)\|^2<\infty,\
  \sum_{|t|=n}\| y(t)\|^2<\infty$, hence  for any $n\in\mathbb{N}\cup
  \{ 0\}$ the series
  $\sum_{|t|=n}[x(t),x(t)]_J$ is absolutely convergent, and
\begin{equation}\label{eq:n-energ}
   \sum\limits_{|t|=n}[x(t),x(t)]_J-\sum\limits_{|t|=n-1}
   [x(t),x(t)]_J=\sum\limits_{|t|=n-1}\| u(t)\|^2-
   \sum\limits_{|t|=n}\| y(t)\|^2;
  \end{equation}
  \item[(ii)] the statement analogous to \textbf{(i)} holds for
  the conjugate system $\alpha^*$.
  \end{description}
\end{prop}
\begin{proof}\textbf{Necessity.} Let the collections
$\{ u(t):|t|\geq 0\}$ and $\{ x(t):|t|=0\}$ of inputs and states
of $\alpha$ satisfy the assumptions of \textbf{(i)}. Apply
induction on $n$. Suppose that for $n=m-1$, where
$m\in\mathbb{N}$, we have $\sum_{|t|=n}\| x(t)\|^2<\infty$. Then,
by virtue of \eqref{eq:n-sys},
\begin{eqnarray*}
\lefteqn{\sum\limits_{|t|=m}\| x(t)\|^2+\sum\limits_{|t|=m}\|
y(t)\|^2=\sum\limits_{|t|=m}\left\| \left(\begin{array}{c}
  x(t) \\
  y(t)
\end{array}\right)\right\|^2} \\
&=& \sum\limits_{|t|=m}\left\|
\sum\limits_{k=1}^NG_k\left(\begin{array}{c}
  x(t-e_k) \\
  u(t-e_k)
\end{array}\right)\right\|^2
\leq   \max_{l\in\{ 1,\ldots,N\} }\|
G_l\|^2\cdot\sum\limits_{|t|=m}\left( \sum\limits_{k=1}^N\left\|
\left(\begin{array}{c}
  x(t-e_k) \\
  u(t-e_k)
\end{array}\right)\right\|\right)^2 \\
&\leq &
 N^2\cdot\max_{l\in\{ 1,\ldots,N\} }\| G_l\|^2\cdot\sum\limits_{|t|=m}
\sum\limits_{k=1}^N\left\|\left(\begin{array}{c}
  x(t-e_k) \\
  u(t-e_k)
\end{array}\right)\right\|^2 \\
&=& N^3\cdot\max_{l\in\{ 1,\ldots,N\} }\| G_l\|^2\cdot
\sum\limits_{|t|=m-1}\left\|\left(\begin{array}{c}
  x(t) \\
  u(t)
\end{array}\right)\right\|^2 \\
&=& N^3\cdot\max_{l\in\{ 1,\ldots,N\} }\| G_l\|^2\left(
\sum\limits_{|t|=m-1}\| x(t)\|^2+\sum\limits_{|t|=m-1}\|
u(t)\|^2\right) <\infty.
\end{eqnarray*}
This implies $\sum_{|t|=n}\| x(t)\|^2<\infty,\
  \sum_{|t|=n}\| y(t)\|^2<\infty$ for $n=m$. Therefore, the latter
  holds for an arbitrary $n\in\mathbb{N}$. Since for any
  $x\in\Hspace{X}$ we have $|[x,x]_J|\leq\| x\|^2$ (see,
   e.g., \cite{AI_1986}),  for any $n\in\mathbb{N}\cup
  \{ 0\}$ the series
  $\sum_{|t|=n}[x(t),x(t)]_J$ converges absolutely.
  Now, by virtue of \eqref{eq:n-sys} and from \eqref{eq:j-c1},
  \eqref{eq:j-c2} we get for any $n\in\mathbb{N}$
\begin{eqnarray*}
\lefteqn{\sum\limits_{|t|=n}[x(t),x(t)]_J+\sum\limits_{|t|=n}\|
y(t)\|^2=\sum\limits_{|t|=n}\left[\left(\begin{array}{c}
  x(t) \\
  y(t)
\end{array}\right),\left(\begin{array}{c}
  x(t) \\
  y(t)
\end{array}\right)\right]_{J_2} } \\
&=&\sum\limits_{|t|=n}\left\langle J_2\left(\begin{array}{c}
  x(t) \\
  y(t)
\end{array}\right),\left(\begin{array}{c}
  x(t) \\
  y(t)
\end{array}\right)\right\rangle \\
 &=& \sum\limits_{|t|=n}\left\langle J_2
\sum\limits_{k=1}^NG_k\left(\begin{array}{c}
  x(t-e_k) \\
  u(t-e_k)
\end{array}\right),\sum\limits_{j=1}^NG_j\left(\begin{array}{c}
  x(t-e_j) \\
  u(t-e_j)
\end{array}\right)\right\rangle \\
&= & \sum\limits_{|t|=n}
\sum\limits_{k=1}^N\sum\limits_{j=1}^N\left\langle G_j^*J_2G_k
\left(\begin{array}{c}
  x(t-e_k) \\
  u(t-e_k)
\end{array}\right),\left(\begin{array}{c}
  x(t-e_j) \\
  u(t-e_j)
\end{array}\right)\right\rangle \\
&=& \sum\limits_{|t|=n}\sum\limits_{k=1}^N\left\langle G_k^*J_2G_k
\left(\begin{array}{c}
  x(t-e_k) \\
  u(t-e_k)
\end{array}\right),\left(\begin{array}{c}
  x(t-e_k) \\
  u(t-e_k)
\end{array}\right)\right\rangle \\
&=& \sum\limits_{|t|=n-1}\left\langle\sum\limits_{k=1}^N
G_k^*J_2G_k \left(\begin{array}{c}
  x(t) \\
  u(t)
\end{array}\right),\left(\begin{array}{c}
  x(t) \\
  u(t)
\end{array}\right)\right\rangle =\sum\limits_{|t|=n-1}\left\langle
J_1\left(\begin{array}{c}
  x(t) \\
  u(t)
\end{array}\right),\left(\begin{array}{c}
  x(t) \\
  u(t)
\end{array}\right)\right\rangle \\
&=& \sum\limits_{|t|=n-1}\left[\left(\begin{array}{c}
  x(t) \\
  u(t)
\end{array}\right),\left(\begin{array}{c}
  x(t) \\
  u(t)
\end{array}\right)\right]_{J_1}=\sum\limits_{|t|=n-1}[x(t),x(t)]_J
+\sum\limits_{|t|=n-1}\| u(t)\|^2,
\end{eqnarray*}
that is equivalent to \eqref{eq:n-energ}, and we have proved the
necessity of condition \textbf{(i)}. The necessity of condition
\textbf{(ii)} is established analogously, by rewriting
\eqref{eq:n-sys} for $\alpha^*$ and using \eqref{eq:j-c3} and
\eqref{eq:j-c4}.

\textbf{Sufficiency.} Let us set for arbitrary $x_0\in\Hspace{X},\
u_0\in\Hspace{U}$
\begin{eqnarray*}
x(t) &:=& \left\{\begin{array}{l}
  x_0 \quad\mbox{for}\ t=0, \\
  0\quad\mbox{anywhere else for}\ |t|=0,
\end{array}\right. \\
u(t) &:=& \left\{\begin{array}{l}
  u_0 \quad\mbox{for}\ t=0, \\
  0\quad\mbox{anywhere else for}\ |t|\geq 0.
\end{array}\right.
\end{eqnarray*}
Clearly, the collections $\{ u(t):|t|\geq 0\}$ and $\{
x(t):|t|=0\}$ of inputs and states of $\alpha$ satisfy the
assumptions of \textbf{(i)}. Then we can write down for them
\eqref{eq:n-energ}, with $n=1$, as follows:
\begin{displaymath}
\sum\limits_{k=1}^N[x(e_k),x(e_k)]_J-[x_0,x_0]_J=\|
u_0\|^2-\sum\limits_{k=1}^N\| y(e_k)\|^2,
\end{displaymath}
or equivalently,
\begin{displaymath}
\left[\left(\begin{array}{c}
  x_0 \\
  u_0
\end{array}\right),\left(\begin{array}{c}
  x_0 \\
  u_0
\end{array}\right)\right]_{J_1}=\sum\limits_{k=1}^N\left[
\left(\begin{array}{c}
  x(e_k) \\
  y(e_k)
\end{array}\right),\left(\begin{array}{c}
  x(e_k) \\
  y(e_k)
\end{array}\right)\right]_{J_2}.
\end{displaymath}
By using system equations \eqref{eq:n-sys}, we obtain

\begin{eqnarray*}
\lefteqn{\left\langle J_1\left(\begin{array}{c}
  x_0 \\
  u_0
\end{array}\right),\left(\begin{array}{c}
  x_0 \\
  u_0
\end{array}\right)\right\rangle =\left[\left(\begin{array}{c}
  x_0 \\
  u_0
\end{array}\right),\left(\begin{array}{c}
  x_0 \\
  u_0
\end{array}\right)\right]_{J_1}=\sum\limits_{k=1}^N\left[
\left(\begin{array}{c}
  x(e_k) \\
  y(e_k)
\end{array}\right),\left(\begin{array}{c}
  x(e_k) \\
  y(e_k)
\end{array}\right)\right]_{J_2} } \\
&=& \sum\limits_{k=1}^N\left\langle J_2 \left(\begin{array}{c}
  x(e_k) \\
  y(e_k)
\end{array}\right),\left(\begin{array}{c}
  x(e_k) \\
  y(e_k)
\end{array}\right)\right\rangle =\sum\limits_{k=1}^N
\left\langle J_2G_k \left(\begin{array}{c}
  x_0 \\
  u_0
\end{array}\right),G_k\left(\begin{array}{c}
  x_0 \\
  u_0
\end{array}\right)\right\rangle \\
&=& \left\langle\sum\limits_{k=1}^N
G_k^*J_2G_k\left(\begin{array}{c}
  x_0 \\
  u_0
\end{array}\right),\left(\begin{array}{c}
  x_0 \\
  u_0
\end{array}\right)\right\rangle .
\end{eqnarray*}
Since $x_0\in\Hspace{X},\ u_0\in\Hspace{U}$ are arbitrary, and the
operators $J_1$ and $\sum_{k=1}^NG_k^*J_2G_k$ are bounded and
selfadjoint, the latter implies \eqref{eq:j-c1}. Analogously,
\textbf{(ii)} implies \eqref{eq:j-c3}.

Now, for arbitrary $x_1,x_2\in\Hspace{X},\ u_1,u_2\in\Hspace{U}$,
and $k,j\in\{ 1,\ldots,N\}\ (k\neq j)$ set
\begin{eqnarray*}
x(t) &:=& \left\{\begin{array}{l}
  x_1 \quad\mbox{for}\ t=e_k-e_j, \\
  x_2 \quad\mbox{for}\ t=0, \\
  0   \quad\mbox{anywhere else for}\ |t|=0,
\end{array}\right. \\
u(t) &:=& \left\{\begin{array}{l} u_1 \quad\mbox{for}\ t=e_k-e_j,
\\
  u_2 \quad\mbox{for}\ t=0, \\
  0\quad\mbox{anywhere else for}\ |t|\geq 0.
\end{array}\right.
\end{eqnarray*}
Clearly, the collections $\{ u(t):|t|\geq 0\}$ and $\{
x(t):|t|=0\}$ of inputs and states of $\alpha$ satisfy the
assumptions of \textbf{(i)}. Then we can write down for them
\eqref{eq:n-energ}, with $n=1$, as follows:
\begin{eqnarray*}
&\sum\limits_{l=1}^N[x(e_k-e_j+e_l),x(e_k-e_j+e_l)]_J+
\sum\limits_{l:\ l\neq k}[x(e_l),x(e_l)]_J-[x_1,x_1]_J-[x_2,x_2]_J
&
\\ &= \| u_1\|^2+\| u_2\|^2-
\sum\limits_{l=1}^N\| y(e_k-e_j+e_l)\|^2- \sum\limits_{l:\ l\neq
k}\| y(e_l)\|^2,&
\end{eqnarray*}
or equivalently,
\begin{eqnarray*}
&\sum\limits_{l=1}^N\left[ \left(\begin{array}{c}
  x(e_k-e_j+e_l) \\
  y(e_k-e_j+e_l)
\end{array}\right),\left(\begin{array}{c}
  x(e_k-e_j+e_l) \\
  y(e_k-e_j+e_l)
\end{array}\right)\right]_{J_2}+ \sum\limits_{l:\ l\neq k}
\left[\left(\begin{array}{c}
  x(e_l) \\
  y(e_l)
\end{array}\right),\left(\begin{array}{c}
  x(e_l) \\
  y(e_l)
\end{array}\right)\right]_{J_2} & \\
 &= \left[\left(\begin{array}{c}
  x_1 \\
  u_1
\end{array}\right),\left(\begin{array}{c}
  x_1 \\
  u_1
\end{array}\right)\right]_{J_1}+\left[\left(\begin{array}{c}
  x_2 \\
  u_2
\end{array}\right),\left(\begin{array}{c}
  x_2 \\
  u_2
\end{array}\right)\right]_{J_1}. &
\end{eqnarray*}
By using system equations \eqref{eq:n-sys}, we obtain
\begin{eqnarray*}
\lefteqn{\sum\limits_{l:\ l\neq j}\left\langle J_2G_l
\left(\begin{array}{c}
  x_1 \\
  u_1
\end{array}\right),G_l\left(\begin{array}{c}
  x_1 \\
  u_1
\end{array}\right)\right\rangle } \\
&+&  \left\langle J_2\left( G_j \left(\begin{array}{c}
  x_1 \\
  u_1
\end{array}\right)+G_k\left(\begin{array}{c}
  x_2 \\
  u_2
\end{array}\right)\right),G_j
\left(\begin{array}{c}
  x_1 \\
  u_1
\end{array}\right)+G_k\left(\begin{array}{c}
  x_2 \\
  u_2
\end{array}\right)\right\rangle  \\
& + & \sum\limits_{l:\ l\neq k} \left\langle
J_2G_l\left(\begin{array}{c}
  x_2 \\
  u_2
\end{array}\right),G_l\left(\begin{array}{c}
  x_2 \\
  u_2
\end{array}\right)\right\rangle=\left\langle J_1\left(\begin{array}{c}
  x_1 \\
  u_1
\end{array}\right),\left(\begin{array}{c}
  x_1 \\
  u_1
\end{array}\right)\right\rangle+\left\langle J_2\left(\begin{array}{c}
  x_2 \\
  u_2
\end{array}\right),\left(\begin{array}{c}
  x_2 \\
  u_2
\end{array}\right)\right\rangle,
\end{eqnarray*}
that is equivalent to
\begin{eqnarray*}
\lefteqn{\left\langle\sum\limits_{l=1}^N G_l^*J_2G_l
\left(\begin{array}{c}
  x_1 \\
  u_1
\end{array}\right),\left(\begin{array}{c}
  x_1 \\
  u_1
\end{array}\right)\right\rangle +\left\langle
\sum\limits_{l=1}^N G_l^*J_2G_l
\left(\begin{array}{c}
  x_2 \\
  u_2
\end{array}\right),\left(\begin{array}{c}
  x_2 \\
  u_2
\end{array}\right)\right\rangle } \\
& + & 2\mbox{Re}\left\langle G_j^*J_2G_k\left(\begin{array}{c}
  x_2 \\
  u_2
\end{array}\right) ,\left(\begin{array}{c}
  x_1 \\
  u_1
\end{array}\right)\right\rangle =\left\langle J_1\left(\begin{array}{c}
  x_1 \\
  u_1
\end{array}\right),\left(\begin{array}{c}
  x_1 \\
  u_1
\end{array}\right)\right\rangle +\left\langle J_1\left(\begin{array}{c}
  x_2 \\
  u_2
\end{array}\right),\left(\begin{array}{c}
  x_2 \\
  u_2
\end{array}\right)\right\rangle .
\end{eqnarray*}
By \eqref{eq:j-c1} established formerly, we obtain:
\begin{displaymath}
2\mbox{Re}\left\langle G_j^*J_2G_k\left(\begin{array}{c}
  x_2 \\
  u_2
\end{array}\right) ,\left(\begin{array}{c}
  x_1 \\
  u_1
\end{array}\right)\right\rangle =0.
\end{displaymath}
One can substitute $-i\tilde{x}_2$ instead of $x_2$, and
$-i\tilde{u}_2$ instead of $u_2$, and obtain
\begin{displaymath}
2\mbox{Im}\left\langle G_j^*J_2G_k\left(\begin{array}{c}
  \tilde{x}_2 \\
  \tilde{u}_2
\end{array}\right) ,\left(\begin{array}{c}
  x_1 \\
  u_1
\end{array}\right)\right\rangle =0.
\end{displaymath}
Since $x_1,x_2,\tilde{x}_2\in\Hspace{X},\
u_1,u_2,\tilde{u}_2\in\Hspace{U}$ can be taken arbitrary,
\eqref{eq:j-c2} follows. Analogously, \textbf{(ii)} implies
\eqref{eq:j-c4}.

The proof is complete.
\end{proof}

Now let us formulate two main theorems of this paper.
\begin{thm}\label{thm:j-c-dil}
An arbitrary system $\alpha =(N;\mathbf{A}, \mathbf{B},
\mathbf{C}, \mathbf{D};\Hspace{X}, \Hspace{U}, \Hspace{Y})$ of the
form \eqref{eq:n-sys} has a dilation $\widetilde{\alpha}
=(N;\widetilde{\mathbf{A}}, \widetilde{\mathbf{B}},
\widetilde{\mathbf{C}}, \mathbf{D};\widetilde{\Hspace{X}},
\Hspace{U}, \Hspace{Y})$, which is a multiparametric
$J$-conservative scattering system for certain canonical symmetry
$J\in L(\widetilde{\Hspace{X}})$.
\end{thm}
\begin{thm}\label{thm:j-c-realiz}
An arbitrary  $L(\Hspace{U}, \Hspace{Y})$-valued function $\theta$
holomorphic on some neighbourhood $\Gamma$ of $z=0$ in
$\nspace{C}{N}$ and vanishing at $z=0$ can be realized as the
transfer function of some system $\dot{\alpha}
=(N;\dot{\mathbf{A}}, \dot{\mathbf{B}}, \dot{\mathbf{C}},
\dot{\mathbf{D}};\dot{\Hspace{X}}, \Hspace{U}, \Hspace{Y})$ of the
form \eqref{eq:n-sys}, which is a multiparametric $J$-conservative
scattering system for certain canonical symmetry $J\in
L(\dot{\Hspace{X}})$, i.e., $\theta(z)=\theta_{\dot{\alpha}}(z)$
in some neighbourhood (possibly, smaller than $\Gamma$) of $z=0$.
\end{thm}

\section{Proofs of the main results}\label{sec:proofs}
In this section we will use the results from
\cite{Ag_1990,K1_2000,K2_2000} on the \emph{Agler--Schur class}
$S_N(\Hspace{U}, \Hspace{Y})$. This class consists of all
$L(\Hspace{U}, \Hspace{Y})$-valued functions
\begin{displaymath}
\theta(z)=\sum\limits_{t\in\nspace{Z}{N}_+}z^t\hat{\theta}_t
\end{displaymath}
(here $\nspace{Z}{N}_+:=\{ t\in\nspace{Z}{N}: t_k\geq 0,\
k=1,\ldots,N\}$, for $z\in\nspace{C}{N}$ and $t\in\nspace{Z}{N}_+$
as usually $z^t:=\prod_{k=1}^Nz_k^{t_k}$, and $\hat{\theta}_t$ are
Taylor coefficients of $\theta$), which are holomorphic on the
open unit polydisk $\nspace{D}{N}:=\{ z\in\nspace{C}{N}:|z_k|<1,\
k=1,\ldots,N\}$ and satisfy the following condition: for any
separable Hilbert space $\Hspace{H}$, any $N$-tuple
$\mathbf{T}=(T_1,\ldots,T_N)$ of commuting contractions on
$\Hspace{H}$, and any positive real $r<1$ one has
\begin{displaymath}
\|\theta(r\mathbf{T})\| <1,
\end{displaymath}
where
\begin{equation}\label{eq:func-calc}
\theta(r\mathbf{T})=\theta(rT_1,\ldots,rT_N):=\sum\limits_
{t\in\nspace{Z}{N}_+}r^{|t|}\mathbf{T}^t\otimes\hat{\theta}_t\in
L(\Hspace{H}\otimes\Hspace{U},\Hspace{H}\otimes\Hspace{Y}),
\end{equation}
$\mathbf{T}^t:=\prod_{k=1}^NT_k^{t_k}$, and the series in
\eqref{eq:func-calc} converges in operator norm.
\begin{lem}\label{lem:lem}
For an arbitrary multiparametric linear system  $\alpha
=(N;\mathbf{A}, \mathbf{B}, \mathbf{C}, \mathbf{D};\Hspace{X},
\Hspace{U}, \Hspace{Y})$ there exist separable Hilbert spaces
$\Hspace{M}_k$, with canonical symmetries $J^{(k)}\in
L(\Hspace{M}_k)$, and holomorphic
$L(\Hspace{X}\oplus\Hspace{U},\Hspace{M}_k)$-valued functions
$F_k$ on $\nspace{D}{N}\ (k=1,\ldots,N)$ such that $
\forall\lambda\in\nspace{D}{N},\ \forall z\in\nspace{D}{N}$
\begin{equation}\label{eq:j-id}
  I_{\Hspace{X}\oplus\Hspace{U}}-(\lambda\mathbf{G})^*(z\mathbf{G})=
  \sum\limits_{k=1}^N(1-\bar{\lambda}_kz_k)F_k(\lambda)^*J^{(k)}
  F_k(z).
\end{equation}
\end{lem}
\begin{proof}
Set $L_{\mathbf{G}}(z):=z\mathbf{G},\
\varepsilon:=\sup_{\mathbf{T}}\|\sum_{k=1}^NT_k\otimes G_k\|$
where this supremum is taken over all $N$-tuples of commuting
contractions $\mathbf{T}=(T_1,\ldots,T_N)$ on a common separable
Hilbert space $\Hspace{H}$. If $\varepsilon\leq 1$ then
$L_{\mathbf{G}}\in
S_N(\Hspace{X}\oplus\Hspace{U},\Hspace{X}\oplus\Hspace{Y})$, and
by Theorem~2.6 of \cite{Ag_1990} the assertion of
Lemma~\ref{lem:lem} follows with $J^{(k)}=I_{\Hspace{M}_k}\
(k=1,\ldots,N)$.

Suppose that $\varepsilon >1$. Then
$\varepsilon^{-1}L_{\mathbf{G}}\in
S_N(\Hspace{X}\oplus\Hspace{U},\Hspace{X}\oplus\Hspace{Y})$. By
Theorem~2.6 of \cite{Ag_1990}, there exist separable Hilbert
spaces $\Hspace{M}_k^+$ and holomorphic
$L(\Hspace{X}\oplus\Hspace{U},\Hspace{M}_k^+)$-valued functions
$H_k^+$ on $\nspace{D}{N}\ (k=1,\ldots,N)$ such that
$\forall\lambda\in\nspace{D}{N},\ \forall z\in\nspace{D}{N}$
\begin{displaymath}
  I_{\Hspace{X}\oplus\Hspace{U}}-(\varepsilon^{-1}\cdot\lambda
  \mathbf{G})^*(\varepsilon^{-1}\cdot z\mathbf{G})=
  \sum\limits_{k=1}^N(1-\bar{\lambda}_kz_k)H_k^+(\lambda)^*H_k^+(z).
\end{displaymath}
Setting $F_k^+(z):=\varepsilon H_k^+(z)$ for $z\in\nspace{D}{N},\
k=1,\ldots,N$, we obtain $\forall\lambda\in\nspace{D}{N},\ \forall
z\in\nspace{D}{N}$:
\begin{equation}\label{eq:id+}
    \varepsilon^2I_{\Hspace{X}\oplus\Hspace{U}}-(\lambda\mathbf{G})^*
  (z\mathbf{G})=
  \sum\limits_{k=1}^N(1-\bar{\lambda}_kz_k)F_k^+(\lambda)^*F_k^+(z).
\end{equation}
Since $\varepsilon^{-1}I_{\Hspace{X}\oplus\Hspace{U}}\in
S_N(\Hspace{X}\oplus\Hspace{U},\Hspace{X}\oplus\Hspace{U})$, again
by Theorem~2.6 of \cite{Ag_1990}, there exist separable Hilbert
spaces $\Hspace{M}_k^-$ and holomorphic
$L(\Hspace{X}\oplus\Hspace{U},\Hspace{M}_k^-)$-valued functions
$H_k^-$ on $\nspace{D}{N}\ (k=1,\ldots,N)$ such that
$\forall\lambda\in\nspace{D}{N},\ \forall z\in\nspace{D}{N}$
\begin{displaymath}
  (1-\varepsilon^{-2})I_{\Hspace{X}\oplus\Hspace{U}}=
  \sum\limits_{k=1}^N(1-\bar{\lambda}_kz_k)H_k^-(\lambda)^*H_k^-(z).
\end{displaymath}
Setting $F_k^-(z):=\varepsilon H_k^-(z)$ for $z\in\nspace{D}{N},\
k=1,\ldots,N$, we obtain $\forall\lambda\in\nspace{D}{N},\ \forall
z\in\nspace{D}{N}$:
\begin{equation}\label{eq:id-}
    (\varepsilon^2-1)I_{\Hspace{X}\oplus\Hspace{U}}=
  \sum\limits_{k=1}^N(1-\bar{\lambda}_kz_k)F_k^-(\lambda)^*F_k^-(z).
\end{equation}

Set $\Hspace{M}_k:=\Hspace{M}_k^+\oplus\Hspace{M}_k^-$, and
according to this orthogonal decomposition define
$F_k:\nspace{D}{N}\to L(\Hspace{X}\oplus\Hspace{U},\Hspace{M}_k)$
by
\begin{displaymath}
F_k(z):=\left(\begin{array}{c}
  F_k^+(z) \\
  F_k^-(z)
\end{array}\right)\quad (z\in\nspace{D}{N}),
\end{displaymath}
and $J^{(k)}:=I_{\Hspace{M}_k^+}\oplus (-I_{\Hspace{M}_k^-})\in
L(\Hspace{M}_k^+\oplus\Hspace{M}_k^-)=L(\Hspace{M}_k)$ for
$k=1,\ldots,N$. By subtracting \eqref{eq:id-} from \eqref{eq:id+}
for each $\lambda\in\nspace{D}{N}$ and $z\in\nspace{D}{N}$, we
obtain \eqref{eq:j-id}, that completes the proof.
\end{proof}

As a by-product of Lemma~\ref{lem:lem}, we obtain the following
result.
\begin{prop}
For an arbitrary $L(\Hspace{U}, \Hspace{Y})$-valued function
$\theta$ holomorphic on some neighbourhood $\Gamma$ of $z=0$ in
$\nspace{C}{N}$ and vanishing at $z=0$ there exist separable
Hilbert spaces $\Hspace{M}_k$, with canonical symmetries
$J^{(k)}\in L(\Hspace{M}_k)$, and holomorphic
$L(\Hspace{U},\Hspace{M}_k)$-valued functions $H_k$ on a
neighbourhood $\Gamma_0\subset\Gamma$ of $z=0$ in $\nspace{C}{N}\
(k=1,\ldots,N)$ such that $ \forall\lambda\in\Gamma_0,\ \forall
z\in\Gamma_0$
\begin{equation}\label{eq:gen-id}
   I_{\Hspace{U}}-\theta(\lambda)^*\theta(z)=
  \sum\limits_{k=1}^N(1-\bar{\lambda}_kz_k)H_k(\lambda)^*J^{(k)}
  H_k(z).
\end{equation}
\end{prop}
\begin{proof}
By Theorem~1 of \cite{K3_2000}, $\theta$ can be realized as the
transfer function of some multiparametric linear system  $\alpha
=(N;\mathbf{A}, \mathbf{B}, \mathbf{C}, \mathbf{D};\Hspace{X},
\Hspace{U}, \Hspace{Y})$, i.e., $\theta(z)=\theta_{\alpha}(z)$ in
some neighbourhood $\Omega\subset\Gamma$ of $z=0$ in
$\nspace{C}{N}$. If $\hat{u}(z)=\sum_{t\in\nspace{Z}{N}_+}z^tu(t)$
is a $\Hspace{U}$-valued function holomorphic on some
neighbourhood of $z=0$ (here $u(\cdot)$ is some input
multisequence of $\alpha$, which has the support in
$\nspace{Z}{N}_+$) then (see \cite{K1_2000}) one can write down
the so-called $Z$-transform of $\alpha$:
\begin{equation}\label{eq:z-tr}
  \hat{\alpha}:\left\{\begin{array}{c}
    \hat{x}(z)=z\mathbf{A}\hat{x}(z)+z\mathbf{B}\hat{u}(z), \\
    \hat{y}(z)=z\mathbf{C}\hat{x}(z)+z\mathbf{D}\hat{u}(z),
  \end{array}\right.
\end{equation}
with holomorphic functions
\begin{eqnarray}
\hat{x}(z) &=&
(I_{\Hspace{X}}-z\mathbf{A})^{-1}z\mathbf{B}\hat{u}(z),
\label{eq:x-tr} \\ \hat{y}(z) &=& \theta_{\alpha}(z)\hat{u}(z)
\label{eq:y-tr}
\end{eqnarray}
on some neighbourhood $\Omega_0$ of $z=0$ in $\nspace{C}{N}$. We
can consider without a loss of generality that
$\Omega_0\subset\Omega$.

Set $\Gamma_0:=\Omega_0\cap\nspace{D}{N}$, and let
$u_1,u_2\in\Hspace{U}$ be arbitrary. Then by using
Lemma~\ref{lem:lem}, and equalities \eqref{eq:z-tr},
\eqref{eq:x-tr}, \eqref{eq:y-tr} twice, for $\hat{u}(\cdot)\equiv
u_1$ and for $\hat{u}(\cdot)\equiv u_2$, we have for all
$\lambda\in\Gamma_0,\ z\in\Gamma_0$:
\begin{eqnarray*}
\lefteqn{\langle\left(I_{\Hspace{U}}-\theta(\lambda)^*\theta(z)
\right)u_1,u_2\rangle_{\Hspace{U}}=\langle
u_1,u_2\rangle_{\Hspace{U}}-\langle\theta(z)u_1,\theta(\lambda)u_2
\rangle_{\Hspace{Y}}} \\ &=& \left\langle\left(\begin{array}{c}
  \hat{x}_1(z) \\
  u_1
\end{array}\right),\left(\begin{array}{c}
  \hat{x}_2(\lambda) \\
  u_2
\end{array}\right)\right\rangle_{\Hspace{X}\oplus\Hspace{U}}-
\left\langle\left(\begin{array}{c}
  \hat{x}_1(z) \\
  \hat{y}_1(z)
\end{array}\right),\left(\begin{array}{c}
  \hat{x}_2(\lambda) \\
  \hat{y}_2(\lambda)
\end{array}\right)\right\rangle_{\Hspace{X}\oplus\Hspace{Y}} \\
&=& \left\langle\left(\begin{array}{c}
  \hat{x}_1(z) \\
  u_1
\end{array}\right),\left(\begin{array}{c}
  \hat{x}_2(\lambda) \\
  u_2
\end{array}\right)\right\rangle_{\Hspace{X}\oplus\Hspace{U}}-
\left\langle z\mathbf{G}\left(\begin{array}{c}
  \hat{x}_1(z) \\
  u_1
\end{array}\right),\lambda\mathbf{G}\left(\begin{array}{c}
  \hat{x}_2(\lambda) \\
  u_2
\end{array}\right)\right\rangle_{\Hspace{X}\oplus\Hspace{Y}} \\
&=& \left\langle\left(
I_{\Hspace{X}\oplus\Hspace{U}}-(\lambda\mathbf{G})^*(z\mathbf{G})
\right)\left(\begin{array}{c}
  \hat{x}_1(z) \\
  u_1
\end{array}\right),\left(\begin{array}{c}
  \hat{x}_2(\lambda) \\
  u_2
\end{array}\right)\right\rangle_{\Hspace{X}\oplus\Hspace{U}} \\
&=& \left\langle\sum\limits_{k=1}^N(1-\bar{\lambda}_kz_k)
F_k(\lambda)^*J^{(k)}F_k(z)\left(\begin{array}{c}
  \hat{x}_1(z) \\
  u_1
\end{array}\right),\left(\begin{array}{c}
  \hat{x}_2(\lambda) \\
  u_2
\end{array}\right)\right\rangle_{\Hspace{X}\oplus\Hspace{U}} \\
&=& \left\langle\sum\limits_{k=1}^N(1-\bar{\lambda}_kz_k)
\left(\begin{array}{c}
  (I_\Hspace{X}-\lambda\mathbf{A})^{-1}\lambda\mathbf{B} \\
  I_\Hspace{U}
\end{array}\right)^*F_k(\lambda)^*J^{(k)}F_k(z)\left(\begin{array}{c}
  (I_\Hspace{X}-z\mathbf{A})^{-1}z\mathbf{B} \\
I_\Hspace{U}
\end{array}\right) u_1, u_2\right\rangle_\Hspace{U}.
\end{eqnarray*}
Setting
\begin{displaymath}
H_k(z):=F_k(z)\left(\begin{array}{c}
   (I_\Hspace{X}-z\mathbf{A})^{-1}z\mathbf{B} \\
I_\Hspace{U}
\end{array}\right)\quad (z\in\Gamma_0),
\end{displaymath}
we obtain \eqref{eq:gen-id}.
\end{proof}
\begin{proof}[Proof of Theorem~\ref{thm:j-c-dil}] Let, as in the
proof of Lemma~\ref{lem:lem}, $L_{\mathbf{G}}(z):=z\mathbf{G},\
\varepsilon:=\sup_{\mathbf{T}}\|\sum_{k=1}^NT_k\otimes G_k\|$
where this supremum is taken over all $N$-tuples of commuting
contractions $\mathbf{T}=(T_1,\ldots,T_N)$ on a common separable
Hilbert space $\Hspace{H}$, and $\mathbf{G}=(G_1,\ldots,G_N)$ be
the $N$-tuple of operators \eqref{eq:gk} corresponding to a given
multiparametric linear system  $\alpha =(N;\mathbf{A}, \mathbf{B},
\mathbf{C}, \mathbf{D};\Hspace{X}, \Hspace{U}, \Hspace{Y})$. If
$\varepsilon\leq 1$ then $L_{\mathbf{G}}\in
S_N(\Hspace{X}\oplus\Hspace{U},\Hspace{X}\oplus\Hspace{Y})$, and
by Theorem~4.2 of \cite{K2_2000}, system $\alpha$ has a
conservative dilation $\widetilde{\alpha}
=(N;\widetilde{\mathbf{A}}, \widetilde{\mathbf{B}},
\widetilde{\mathbf{C}}, \mathbf{D};\widetilde{\Hspace{X}},
\Hspace{U}, \Hspace{Y})$, i.e. a $J$-conservative one with
$J=I_{\widetilde{\Hspace{X}}}$.

Suppose now that $\varepsilon >1$. Applying Lemma~\ref{lem:lem} to
$\alpha$, we have the existence of separable Hilbert spaces
$\Hspace{M}_k$ with canonical symmetries $J^{(k)}\in
L(\Hspace{M}_k)$, and holomorphic
$L(\Hspace{X}\oplus\Hspace{U},\Hspace{M}_k)$-valued functions
$F_k$ on $\nspace{D}{N}\ (k=1,\ldots,N)$ such that \eqref{eq:j-id}
holds. Let us define these spaces $\Hspace{M}_k$, operators
$J^{(k)}$, and functions $F_k\ (k=1,\ldots,N)$ exactly as in the
proof of Lemma~\ref{lem:lem}, i.e.,
$\Hspace{M}_k:=\Hspace{M}_k^+\oplus\Hspace{M}_k^-,\
J^{(k)}:=I_{\Hspace{M}_k^+}\oplus (-I_{\Hspace{M}_k^-})\in
L(\Hspace{M}_k^+\oplus\Hspace{M}_k^-)=L(\Hspace{M}_k)$,
\begin{displaymath}
F_k(z):=\left(\begin{array}{c}
  F_k^+(z) \\
  F_k^-(z)
\end{array}\right)\in L(\Hspace{X}\oplus\Hspace{U},
\Hspace{M}_k^+\oplus\Hspace{M}_k^-)\quad (z\in\nspace{D}{N}),
\end{displaymath}
 so that
\eqref{eq:id+} and \eqref{eq:id-} hold. Set $\Hspace{M}^\pm
:=\bigoplus_{k=1}^N\Hspace{M}_k^\pm,\
\Hspace{M}:=\bigoplus_{k=1}^N\Hspace{M}_k=
\Hspace{M}^+\oplus\Hspace{M}^-,\
J_\Hspace{M}:=\bigoplus_{k=1}^NJ^{(k)}\in
L(\bigoplus_{k=1}^N\Hspace{M}_k)=L(\Hspace{M})$,
\begin{displaymath}
F^\pm (z):=\left(\begin{array}{c}
  F_1^\pm (z) \\
  \vdots \\
  F_N^\pm (z)
\end{array}\right)\in L(\Hspace{X}\oplus\Hspace{U},
\Hspace{M}^\pm),\quad F(z):=\left(\begin{array}{c}
  F_1(z) \\
  \vdots \\
  F_N(z)
\end{array}\right)\in L(\Hspace{X}\oplus\Hspace{U},
\Hspace{M})\quad (z\in\nspace{D}{N}),
\end{displaymath}
$P_k^\pm :=P_{\Hspace{M}_k^\pm}\in L(\Hspace{M}^\pm ),\ P_k
:=P_{\Hspace{M}_k}\in L(\Hspace{M})\ (k=1,\ldots,N)$. Then
$\forall z\in\nspace{C}{N}\ z\mathbf{P}^\pm
=\bigoplus_{k=1}^Nz_kI_{\Hspace{M}_k^\pm}\in L(\Hspace{M}^\pm ),\
 z\mathbf{P}
=\bigoplus_{k=1}^Nz_kI_{\Hspace{M}_k}\in L(\Hspace{M})$.

It follows from \eqref{eq:id+} that
$\forall\lambda\in\nspace{D}{N},  \forall z\in\nspace{D}{N}$
\begin{equation}\label{eq:f+0}
   F^+(0)^*F^+(0)=F^+(0)^*F^+(z)=F^+(\lambda)^*F^+(0)=\varepsilon^2
  I_{\Hspace{X}\oplus\Hspace{U}}.
\end{equation}
In particular, $F^+(0)$ is a bounded and boundedly invertible
operator, and $F^+(0)(\Hspace{X}\oplus\Hspace{U})$ is a closed
lineal, i.e. a subspace in $\Hspace{M}^+$. Analogously, it follows
from \eqref{eq:id-} that $\forall\lambda\in\nspace{D}{N}, \forall
z\in\nspace{D}{N}$
\begin{equation}\label{eq:f-0}
   F^-(0)^*F^-(0)=F^-(0)^*F^-(z)=F^-(\lambda)^*F^-(0)=
  (\varepsilon^2-1)
  I_{\Hspace{X}\oplus\Hspace{U}}.
\end{equation}
In particular, $F^-(0)$ is a bounded and boundedly invertible
operator, and $F^-(0)(\Hspace{X}\oplus\Hspace{U})$ is a closed
lineal, i.e. a subspace in $\Hspace{M}^-$. It follows from
\eqref{eq:f+0} and \eqref{eq:f-0} that
$\forall\lambda\in\nspace{D}{N},  \forall z\in\nspace{D}{N}$
\begin{equation}\label{eq:f+-0}
   (F^\pm(\lambda)-F^\pm(0))^*(F^\pm(z)-F^\pm(0))=
  F^\pm(\lambda)^*F^\pm(z)-F^\pm(0)^*F^\pm(0).
\end{equation}
Taking into account \eqref{eq:f+0}, rewrite \eqref{eq:id+} as
\begin{displaymath}
F^+(0)^*F^+(0)-(\lambda\mathbf{G})^*(z\mathbf{G})=
  F^+(\lambda)^*F^+(z)-(\lambda\mathbf{P}^+F^+(\lambda))^*
  (z\mathbf{P}^+F^+(z)) \quad (\lambda\in\nspace{D}{N},\
   z\in\nspace{D}{N}),
\end{displaymath}
and by virtue of \eqref{eq:f+-0}, we get
$\forall\lambda\in\nspace{D}{N},  \forall z\in\nspace{D}{N}$:
\begin{equation}\label{eq:id++}
 (\lambda\mathbf{P}^+F^+(\lambda))^*
  (z\mathbf{P}^+F^+(z))=\left(\begin{array}{c}
   F^+(\lambda)-F^+(0)  \\
\lambda\mathbf{G}
  \end{array}\right)^*\left(\begin{array}{c}
   F^+(z)-F^+(0)  \\
z\mathbf{G}
  \end{array}\right).
  \end{equation}
Taking into account \eqref{eq:f-0}, rewrite \eqref{eq:id-} as
\begin{displaymath}
F^-(0)^*F^-(0)=
  F^-(\lambda)^*F^-(z)-(\lambda\mathbf{P}^-F^-(\lambda))^*
  (z\mathbf{P}^-F^-(z)) \quad (\lambda\in\nspace{D}{N},\
   z\in\nspace{D}{N}),
\end{displaymath}
and by virtue of \eqref{eq:f+-0}, we get
$\forall\lambda\in\nspace{D}{N},  \forall z\in\nspace{D}{N}$:
\begin{equation}\label{eq:id--}
 (\lambda\mathbf{P}^-F^-(\lambda))^*
  (z\mathbf{P}^-F^-(z))=(F^-(\lambda)-F^-(0))^*(F^-(z)-F^-(0)).
  \end{equation}
Now, adding \eqref{eq:id--} to \eqref{eq:id++}, we get
$\forall\lambda\in\nspace{D}{N},  \forall z\in\nspace{D}{N}$:
\begin{equation}\label{eq:id}
 (\lambda\mathbf{P}F(\lambda))^*
  (z\mathbf{P}F(z))=\left(\begin{array}{c}
   F(\lambda)-F(0)  \\
\lambda\mathbf{G}
  \end{array}\right)^*\left(\begin{array}{c}
   F(z)-F(0)  \\
z\mathbf{G}
  \end{array}\right),
  \end{equation}
and subtracting \eqref{eq:id--} from \eqref{eq:id++}, we get
$\forall\lambda\in\nspace{D}{N},  \forall z\in\nspace{D}{N}$:
\begin{equation}\label{eq:id-j}
 (\lambda\mathbf{P}F(\lambda))^*J_\Hspace{M}
  (z\mathbf{P}F(z))=\left(\begin{array}{c}
   F(\lambda)-F(0)  \\
\lambda\mathbf{G}
  \end{array}\right)^*\left(\begin{array}{cc}
    J_\Hspace{M} & 0 \\
    0 & I_{\Hspace{X}\oplus\Hspace{Y}} \
  \end{array}\right)\left(\begin{array}{c}
   F(z)-F(0)  \\
z\mathbf{G}
  \end{array}\right).
  \end{equation}

It follows from \eqref{eq:id} that there exists unique unitary
operator
\begin{equation}\label{eq:u1}
  U:\bigvee_{z\in\nspace{D}{N}}z\mathbf{P}F(z)
  (\Hspace{X}\oplus\Hspace{U})\to\bigvee_{z\in\nspace{D}{N}}
  \left(\begin{array}{c}
   F(z)-F(0)  \\
z\mathbf{G}
  \end{array}\right)(\Hspace{X}\oplus\Hspace{U})
\end{equation}
such that $\forall z\in\nspace{D}{N}$:
\begin{equation}\label{eq:u2}
  U(z\mathbf{P})F(z)=\left(\begin{array}{c}
   F(z)-F(0)  \\
z\mathbf{G}
  \end{array}\right).
\end{equation}

It follows from \eqref{eq:f+0} and \eqref{eq:f-0}  that
\begin{equation}\label{eq:semiun}
  F(0)^*J_\Hspace{M}F(0)=I_{\Hspace{X}\oplus\Hspace{U}},
\end{equation}
i.e. $F(0)\in L(\Hspace{X}\oplus\Hspace{U},\Hspace{M})$ is a
$(I_{\Hspace{X}\oplus\Hspace{U}},J_\Hspace{M})$-semiunitary
operator. Moreover,
$F(0)(\Hspace{X}\oplus\Hspace{U})=F^+(0)(\Hspace{X}\oplus\Hspace{U})
\oplus F^-(0)(\Hspace{X}\oplus\Hspace{U})$ is a closed lineal,
i.e. a subspace in $\Hspace{M}=\Hspace{M}^+\oplus\Hspace{M}^-$. In
addition, from \eqref{eq:f+0} and \eqref{eq:f-0} one can see that
$\forall z\in\nspace{D}{N}$:
\begin{displaymath}
  F(0)^*(F(z)-F(0))=0,
\end{displaymath}
hence
\begin{displaymath}
  (F(z)-F(0))(\Hspace{X}\oplus\Hspace{U})\subset\Hspace{M}
  \ominus F(0)(\Hspace{X}\oplus\Hspace{U}).
  \end{displaymath}

Now let us show that the subspace $\Hspace{K}_0:=\Hspace{M}
  \ominus F(0)(\Hspace{X}\oplus\Hspace{U})$ in $\Hspace{M}$ is a
  Krein space with respect to the metric $[\cdot,\cdot]_{J_0}$
  induced by the canonical symmetry $J_0:=P_{\Hspace{K}_0}
  J_\Hspace{M}|\Hspace{K}_0$. By Theorem~I.7.16 of \cite{AI_1986},
  in order that $J_0$ is a canonical symmetry on $\Hspace{K}_0$,
  it is necessary and sufficient that any $h\in\Hspace{M}$ has a
  $J_\Hspace{M}$-orthogonal projection onto $\Hspace{K}_0$, i.e. a
  vector $h_0\in\Hspace{K}_0$ such that $J_\Hspace{M}(h-h_0)\bot
  \Hspace{K}_0$. For an arbitrary $h\in\Hspace{M}$ set $h_0:=h-
  J_\Hspace{M}F(0)F(0)^*h$. Since, due to \eqref{eq:semiun},
\begin{displaymath}
  F(0)^*h_0=F(0)^*h-F(0)^*J_\Hspace{M}F(0)F(0)^*h=0,
  \end{displaymath}
we get $h_0\in\Hspace{M}\ominus
F(0)(\Hspace{X}\oplus\Hspace{U})=\Hspace{K}_0$. For an arbitrary
$g\in\Hspace{K}_0$ we have:
\begin{displaymath}
  \langle J_\Hspace{M}(h-h_0),g\rangle =\langle J_\Hspace{M}^2
  F(0)F(0)^*h,g\rangle =\langle F(0)F(0)^*h,g\rangle=0.
  \end{displaymath}
  Thus, $h_0$ is a desired
  $J_\Hspace{M}$-orthogonal projection of $h$ onto $\Hspace{K}_0$,
  and we have proved that $J_0$ is a canonical symmetry on
   $\Hspace{K}_0$ (i.e., $\Hspace{K}_0$ is a
  Krein space with respect to the metric $[\cdot,\cdot]_{J_0}$).

  Further,
  \begin{eqnarray*}
\bigvee_{z\in\nspace{D}{N}}z\mathbf{P}F(z)
  (\Hspace{X}\oplus\Hspace{U}) &\subset &\Hspace{M},\\
\bigvee_{z\in\nspace{D}{N}}
  \left(\begin{array}{c}
   F(z)-F(0)  \\
z\mathbf{G}
  \end{array}\right)(\Hspace{X}\oplus\Hspace{U}) &\subset &
  \Hspace{K}_0\oplus\Hspace{X}\oplus\Hspace{Y}.
\end{eqnarray*}
Define on the space
$\Hspace{K}_I:=\Hspace{K}_0\oplus\Hspace{X}\oplus\Hspace{Y}$ the
canonical symmetry $J_I:=J_0\oplus
I_{\Hspace{X}\oplus\Hspace{Y}}\in L(\Hspace{K}_I)$. Due to
\eqref{eq:id-j}, the operator $U$, defined by  \eqref{eq:u1},
\eqref{eq:u2}, can be considered as a bounded and boundedly
invertible on its domain $(J_\Hspace{M},J_I)$-isometric operator
from $\Hspace{M}$ to $\Hspace{K}_I$, whose domain and range are
given by
\begin{displaymath}
  \mathfrak{D}(U)=\bigvee_{z\in\nspace{D}{N}}z\mathbf{P}F(z)
  (\Hspace{X}\oplus\Hspace{U}),\quad \mathfrak{R}(U)=
  \bigvee_{z\in\nspace{D}{N}}
  \left(\begin{array}{c}
   F(z)-F(0)  \\
z\mathbf{G}
  \end{array}\right)(\Hspace{X}\oplus\Hspace{U}).
  \end{displaymath}
By  Theorem V.2.18 of \cite{AI_1986}, there exist a separable
Hilbert space $\Hspace{K}_{II}$, a canonical symmetry $J_{II}\in
L(\Hspace{K}_{II})$, and a $(\check{J}_1,\check{J}_2)$-unitary
operator
$\check{U}:\Hspace{K}_{II}\oplus\Hspace{M}\to\Hspace{K}_{II}
\oplus\Hspace{K}_I$, with
\begin{displaymath}
\check{J}_1:=J_{II}\oplus J_\Hspace{M}\in
L(\Hspace{K}_{II}\oplus\Hspace{M}),\quad \check{J}_2:=J_{II}\oplus
J_I\in L(\Hspace{K}_{II}\oplus\Hspace{K}_I)
  \end{displaymath}
such that $\check{U}$ is an extension of $U$, i.e.
$P_{\mathfrak{R}(U)}\check{U}| \mathfrak{D}(U)=U$.

Since for any $z\in\nspace{D}{N}$ we have
\begin{displaymath}
  (F(z)-F(0))(\Hspace{X}\oplus\Hspace{U})\subset\Hspace{M}
  \ominus F(0)(\Hspace{X}\oplus\Hspace{U})=\Hspace{K}_0,
  \end{displaymath}
we get for any $z\in\nspace{D}{N}$:
\begin{displaymath}
  F(z)=\left(\begin{array}{c}
    F(z)-F(0) \\
    F(0)
  \end{array}\right)\in L(\Hspace{X}\oplus\Hspace{U},\Hspace{M})=
  L(\Hspace{X}\oplus\Hspace{U},\Hspace{K}_0\oplus
  F(0)(\Hspace{X}\oplus\Hspace{U})).
  \end{displaymath}

Set for all $k\in\{ 1,\ldots,N\}\quad
\check{P}_k:=\delta_{1k}I_{\Hspace{K}_{II}}\oplus P_k\in
L(\Hspace{K}_{II}\oplus\Hspace{M})$, where $\delta_{1k}$ denotes
the Kronecker symbol,
\begin{eqnarray*}
  \check{G}_k &:=& \check{U}\check{P}_k
  (I_{\Hspace{K}_{II}\oplus\Hspace{K}_0}\oplus F(0))\in
L(\Hspace{K}_{II}\oplus\Hspace{K}_0\oplus\Hspace{X}\oplus\Hspace{U},
\Hspace{K}_{II}\oplus\Hspace{K}_I)\\ &=&
L(\Hspace{K}_{II}\oplus\Hspace{K}_0\oplus\Hspace{X}\oplus\Hspace{U},
\Hspace{K}_{II}\oplus\Hspace{K}_0\oplus\Hspace{X}\oplus\Hspace{Y}).
  \end{eqnarray*}
Clearly, $I_{\Hspace{K}_{II}\oplus\Hspace{K}_0}\oplus F(0)\in
L(\Hspace{K}_{II}\oplus\Hspace{K}_0\oplus\Hspace{X}\oplus\Hspace{U},
\Hspace{K}_{II}\oplus\Hspace{M})$ is a $(J_1,\check{J}_1)$-unitary
operator, where
\begin{displaymath}
 J_1:=J_{II}\oplus J_0\oplus I_{\Hspace{X}\oplus\Hspace{U}}\in
 L(\Hspace{K}_{II}\oplus\Hspace{K}_0\oplus\Hspace{X}\oplus
 \Hspace{U}),
  \end{displaymath}
for any $\zeta\in\nspace{T}{N}\quad \zeta\check{\mathbf{P}}\in
L(\Hspace{K}_{II}\oplus\Hspace{M})$ is a $\check{J}_1$-unitary
operator, and $\check{U}\in L(\Hspace{K}_{II}\oplus\Hspace{M},
\Hspace{K}_{II}\oplus\Hspace{K}_I)$ is a
$(\check{J}_1,J_2)$-unitary operator, where
\begin{displaymath}
 J_2:=\check{J}_2= J_{II}\oplus J_0\oplus I_{\Hspace{X}\oplus
 \Hspace{Y}}\in
 L(\Hspace{K}_{II}\oplus\Hspace{K}_I)=
 L(\Hspace{K}_{II}\oplus\Hspace{K}_0\oplus\Hspace{X}\oplus
 \Hspace{Y}).
  \end{displaymath}
Therefore, for any $\zeta\in\nspace{T}{N}$ the operator
$\zeta\mathbf{G}$ is $(J_1,J_2)$-unitary.

Consider the following partitioning of $\check{G}_k\
(k=1,\ldots,N)$:
\begin{displaymath}
\check{G}_k=\left(\begin{array}{cc}
  \check{A}_k & \check{B}_k \\
  \check{C}_k & \check{D}_k
\end{array}\right) :(\Hspace{K}_{II}\oplus\Hspace{K}_0)\oplus
(\Hspace{X}\oplus \Hspace{U})\to
(\Hspace{K}_{II}\oplus\Hspace{K}_0)\oplus (\Hspace{X}\oplus
\Hspace{Y}).
  \end{displaymath}
Then $\check{\alpha}:=(N;\check{\mathbf{A}},\check{\mathbf{B}},
\check{\mathbf{C}},\check{\mathbf{D}};
\Hspace{K}_{II}\oplus\Hspace{K}_0,\Hspace{X}\oplus
\Hspace{U},\Hspace{X}\oplus \Hspace{Y})$ is a multiparametric
$(J_{II}\oplus J_0)$-conservative scattering system.

 Since, by virtue of \eqref{eq:u2}, $\forall z\in\nspace{D}{N}$
\begin{eqnarray*}
  z\check{\mathbf{G}}\left(\begin{array}{c}
    F(z)-F(0) \\
    I_{\Hspace{X}\oplus\Hspace{U}}
  \end{array}\right) &=& \check{U}(z\check{\mathbf{P}})
  \left(\begin{array}{cc}
    I_{\Hspace{K}_{II}\oplus\Hspace{K}_0} & 0 \\
    0 & F(0)
  \end{array}\right)\left(\begin{array}{c}
    F(z)-F(0) \\
    I_{\Hspace{X}\oplus\Hspace{U}}
  \end{array}\right) \\
&=&
\check{U}(z\check{\mathbf{P}})F(z)=U(z\mathbf{P})F(z)=
\left(\begin{array}{c}
    F(z)-F(0) \\
    z\mathbf{G}
  \end{array}\right),
    \end{eqnarray*}
we have for all $z\in\nspace{D}{N}$:
\begin{eqnarray*}
  z\check{\mathbf{A}}(F(z)-F(0))+z\check{\mathbf{B}} &=&
    F(z)-F(0), \\
z\check{\mathbf{C}}(F(z)-F(0))+z\check{\mathbf{D}} &=&
    z\mathbf{G}.
    \end{eqnarray*}
Therefore, in some neighbourhood $\Omega\subset\nspace{D}{N}$ of
$z=0$ the resolvent
$(I_{\Hspace{K}_{II}\oplus\Hspace{K}_0}-z\check{\mathbf{A}})^{-1}$
is well-defined and holomorphic, and we have in this
neighbourhood:
\begin{eqnarray}
  F(z)-F(0) &=& (I_{\Hspace{K}_{II}\oplus\Hspace{K}_0}-
  z\check{\mathbf{A}})^{-1}z\check{\mathbf{B}}, \nonumber \\
  z\check{\mathbf{G}} &=& z\check{\mathbf{D}}+z\check{\mathbf{C}}
(I_{\Hspace{K}_{II}\oplus\Hspace{K}_0}-
  z\check{\mathbf{A}})^{-1}z\check{\mathbf{B}}=
  \theta_{\check{\alpha}}(z). \label{eq:lin-tf}
    \end{eqnarray}
    Thus, $\check{\alpha}$ is a
    $(J_{II}\oplus J_0)$-conservative scattering system
realization of the linear operator-valued function $L_\mathbf{G}$.
>From \eqref{eq:lin-tf} we obtain:
\begin{eqnarray*}
 \forall z\in\Omega_0\quad z\check{\mathbf{D}} &=&
    z\mathbf{G}\\
 \forall z\in\Omega_0,\ \forall n\in\mathbb{N}\cup\{ 0\}\quad
  z\check{\mathbf{C}}(z\check{\mathbf{A}})^nz\check{\mathbf{B}}
&=& 0.
\end{eqnarray*}
As it was shown in the proof of Theorem~4.2 of \cite{K2_2000} (see
also Subsection~V.3.1 of \cite{AI_1986}), the latter means that
the system $\widetilde{\alpha}:=(N;\widetilde{\mathbf{A}},
\widetilde{\mathbf{B}},\widetilde{\mathbf{C}},\mathbf{D};
\widetilde{\Hspace{X}}=\Hspace{K}_{II}\oplus\Hspace{K}_0\oplus
\Hspace{X},\Hspace{U},\Hspace{Y})$ which is determined by the
system operators $\widetilde{G}_k$ coinciding with the system
operators $\check{G}_k$ of $\check{\alpha}\ (k=1,\ldots,N)$,
however with another block partitioning of these operators, is a
dilation of $\alpha
=(N;\mathbf{A},\mathbf{B},\mathbf{C},\mathbf{D};
\Hspace{X},\Hspace{U},\Hspace{Y})$. Defining a canonical symmetry
of $\widetilde{\Hspace{X}}$:
\begin{displaymath}
J:=J_{II}\oplus J_0\oplus I_{\Hspace{X}}\in
L(\Hspace{K}_{II}\oplus\Hspace{K}_0\oplus\Hspace{X})=
L(\widetilde{\Hspace{X}}),
  \end{displaymath}
  we obtain that $\widetilde{\alpha}$ is a desired
  $J$-conservative scattering system dilation of a given system
  $\alpha$.

  The proof is complete.
\end{proof}

\begin{proof}[Proof of
Theorem~\ref{thm:j-c-realiz}] Let $\theta$ be a given
$L(\Hspace{U},\Hspace{Y})$-valued function holomorphic on a
neighbourhood $\Gamma$ of $z=0$ in $\nspace{C}{N}$ and vanishing
at $z=0$. Then, by Theorem~1 of \cite{K3_2000}, there exists a
realization $\alpha
=(N;\mathbf{A},\mathbf{B},\mathbf{C},\mathbf{D};
\Hspace{X},\Hspace{U},\Hspace{Y})$ of this function, i.e.,
$\theta(z)=\theta_\alpha (z)$ in some neighbourhood
$\Gamma_0\subset\Gamma$ of $z=0$. By Theorem~\ref{thm:j-c-dil} of
this paper, there exists a multiparametric  $J$-conservative
scattering system $\widetilde{\alpha}:=(N;\widetilde{\mathbf{A}},
\widetilde{\mathbf{B}},\widetilde{\mathbf{C}},\mathbf{D};
\widetilde{\Hspace{X}},\Hspace{U},\Hspace{Y})$, with a canonical
symmetry $J\in L(\widetilde{\Hspace{X}})$, which is a dilation of
$\alpha$. By Proposition~3.8 of \cite{K2_2000}, the transfer
functions of a multiparametric linear system of the form
\eqref{eq:n-sys} and of its dilation coincide in some
neighbourhood of $z=0$ in $\nspace{C}{N}$. Hence,
$\theta_{\widetilde{\alpha}}(z)=\theta_\alpha (z)=\theta (z)$ in
some neighbourhood $\Omega_0\subset\Gamma_0$ of $z=0$. Thus, the
system $\dot{\alpha}:=\widetilde{\alpha}$ is a desired
 $J$-conservative scattering system realization of $\theta$.

 The proof is complete.
\end{proof}

 \vspace{1cm} \noindent\textbf{Acknowledgement.} I am
thankful to D.Z.~Arov for suggesting this problem.

    \ \\
Department of Higher Mathematics \\ Odessa State Academy of Civil
Engineering and Architecture \\ Didrihson str. 4, Odessa, 65029,
Ukraine \\
\\
2000 Mathematics Subject Classification: 47A20, 47A48, 47B50,
93C35

\end{document}

%% file: mydef.tex
\def\ifundefined#1{\expandafter\ifx\csname#1\endcsname\relax}
\providecommand{\comment}[1]{}

\usepackage{amsfonts}
\providecommand{\tthdump}[1]{#1}
\tthdump{}
\newcommand{\Cliff}[2][\comment]{{\ensuremath{%
\mathcal{C}\kern-0.18em\ell(#1,#2)}}}

\ifundefined{qed}
    \DeclareMathSymbol{\qed}{0}{AMSa}{"03}
\fi

\providecommand{\eqref}[1]{\textup{(\ref{#1})}}

\providecommand{\href}[2]{#2}
\iftth

\else

\fi